\documentclass{amsart}
\usepackage{lattice}

\DeclareMathOperator{\Con}{Con}
\DeclareMathOperator{\Conc}{Con_c}

\newcommand{\RR}{\mathbb{R}}
\newcommand{\RRh}{\widehat{\mathbb{R}}}

\newcommand{\id}{\operatorname{id}}
\newcommand{\compose}{\circ}
\newcommand{\two}{\boldsymbol{2}}
\newcommand{\famm}[2]{(\,#1\mid#2\,)}
\newcommand{\jz}{$\set{\jj,0}$}
\newcommand{\res}{\!\restriction}
\newcommand{\ot}{\set{1,2}}

\DeclareMathOperator{\J}{J}
\DeclareMathOperator{\M}{M}

\theoremstyle{plain}
\newtheorem{lemma}{Lemma}[section]
\newtheorem{theorem}{Theorem}
\newtheorem{proposition}[lemma]{Proposition}

\newtheorem*{tuma}{T\r uma's Theorem}
\newtheorem*{Huhn}{Huhn's Theorem}
\newtheorem*{Dilworth}{Dilworth' Theorem}
\newtheorem*{tisch}{Tischendorf's Theorem}

\newtheorem*{add}{Addition for $D$ Boolean}

\theoremstyle{definition}

\newtheorem{problem}{Problem}

\numberwithin{equation}{section}

\begin{document}

\title{Congruence amalgamation of lattices}

  \author[G. Gr\"atzer]{George~Gr\"{a}tzer}
  \author[H. Lakser]{Harry~Lakser}

  \address{Department of Mathematics\\
          University of Manitoba\\
          Winnipeg, MB, R3T~2N2\\
          Canada}
  \email[G. Gr\"atzer]{gratzer@cc.umanitoba.ca}
  \urladdr[G. Gr\"atzer]{http://www.maths.umanitoba.ca/homepages/gratzer/}
  \email[H. Lakser]{hlakser@cc.umanitoba.ca}

  \thanks{The research of the first two authors was supported by the NSERC
  of Canada}

  \author[F.~Wehrung]{Friedrich Wehrung}
\address{CNRS, ESA 6081\\
Universit\'e de Caen, Campus II\\
D\'epartement de Math\'ematiques\\
B.P. 5186\\
14032 Caen Cedex\\
France}
 \email[F.~Wehrung]{wehrung@math.unicaen.fr}
 \urladdr[F.~Wehrung]{http://www.math.unicaen.fr/\~{}wehrung}
 \date{October 18, 1999}
 \subjclass{Primary:06B05, 06B10, Secondary:06D05}
 \keywords{Congruence, amalgamation, lattice, distributive, relatively
complemented} 
\begin{abstract}
J.~T\r uma proved an interesting ``congruence amalgamation'' result.  We are
generalizing and providing an alternate proof for it.
We then provide applications of this result:
\begin{enumerate}
\item A.P. Huhn proved that every distributive algebraic lattice $D$ with at
most $\aleph_1$ compact elem\-ents can be represented as the congruence lattice
of a lattice~$L$. We show that $L$ can be constructed as a locally finite
relatively complemented lattice with zero.

\item We find a large class of lattices, the 
\emph{$\go$-congruence-finite} lattices,
that contains all locally finite countable lattices, in which every 
lattice has a
relatively complemented congruence-pre\-serv\-ing extension.

\end{enumerate}
  \end{abstract}
\maketitle

\section{Introduction} \label{S:intro}
The first congruence lattice characterization theorem is due to  R.P. Dilworth
(see G. Gr\"atzer and E.T. Schmidt~\cite{GS62}):

  \begin{Dilworth}
Let $D$ be a finite distributive lattice.  Then there exists a finite
lattice $L$ such that the congruence lattice of $L$, $\Con L$, is isomorphic
to~$D$.
  \end{Dilworth}

The best extension of this result is due to
A.P. Huhn \cite{Huhn89b}:

\begin{Huhn}
Let $D$ be a distributive algebraic lattice. If $D$ has at most
$\aleph_1$ compact elements, then there exists a lattice $L$ such that
$\Con L \iso D$.
  \end{Huhn}

An equivalent form of this result is the following: \emph{Let $S$ be 
a distributive
join-semilattice with zero. If $|S| \leq \aleph_1$, then there exists 
a lattice $L$
such that the join-semilattice of compact congruences of $L$ is 
isomorphic to $S$.}

By P. Pudl\'ak \cite{Pudl85}, $S$ is a direct limit of its finite distributive
\jz-subsemi\-lattices. So it is natural to attempt to prove Huhn's result
with a direct limit argument.

Assigning to a lattice $L$ its congruence lattice, $\Con{L}$, determines a
functor $\Con{}$ from the category of lattices with lattice homomorphisms
to the category of algebraic distributive lattices with morphisms the
complete $\JJ$-homomorphisms. Specifically, if $K$ and $L$ are lattices and
$\gf\colon K\to L$ is a lattice homomorphism, then the mapping
$\Con{\gf}\colon \Con{K}\to \Con{L}$
is determined by setting
\begin{equation*}
(\Con{\gf})\gQ = \gQ_{L}\famm{\vv<\gf x,\gf y>}{x,y\in K,\ 
{x}\equiv{y}\,(\gQ)},
\end{equation*}
for each $\gQ\in\Con{K}$.

J.~T\r uma \cite{Tuma} proved the following result:

\begin{tuma}
Let $L_0$, $L_1$, $L_2$ be finite atomistic lattices and let $\gh_1\colon
L_0\to L_1$ and $\gh_2\colon L_0\to L_2$ be lattice
embeddings preserving the zero
such that $\Con{\gh_1}$ and $\Con{\gh_2}$ are injective.
Let $D$ be a finite distributive lattice, and, for $i\in\ot$, let
$\gy_i\colon\Con{L_i}\to D$ be \jz-embeddings such that
\begin{equation*}
\gy_1\compose\Con{\gh_1} = \gy_2\compose\Con{\gh_2}.
\end{equation*}
Then there is a finite atomistic lattice $L$, and there are
lattice embeddings
$\gf_i\colon L_i\to L$, for $i\in\ot$, that preserve the zero, satisfying
\begin{equation*}
\gf_1\compose\gh_1=\gf_2\compose\gh_2,
\end{equation*}
and there is an isomorphism $\ga\colon\Con{L}\to D$ such that
\begin{equation*}
\ga\compose\Con{\gf_i} = \gy_i, \q\text{ for } i\in\ot.
\end{equation*}
\end{tuma}

We extend T\r uma's result by proving:

\begin{theorem}\label{T:Tumagen}
Let $L_0$, $L_1$, $L_2$ be lattices and let $\gh_1\colon L_0\to L_1$ and
$\gh_2\colon L_0\to L_2$ be lattice homomorphisms. Let $D$ be a finite
distributive lattice, and, for $i\in\ot$, let $\gy_i\colon \Con{L_i}\to D$
be complete $\JJ$-homomorphisms such that
\begin{equation*}
\gy_1\compose\Con{\gh_1} = \gy_2\compose\Con{\gh_2}.
\end{equation*}
There is then a lattice $L$, there are lattice homomorphisms
$\gf_i\colon L_i\to L$, for $i\in\ot$, with
\begin{equation*}
\gf_1\compose\gh_1=\gf_2\compose\gh_2,
\end{equation*}
and there is an isomorphism $\ga\colon\Con{L}\to D$ such that
\begin{equation*}
\ga\compose\Con{\gf_i} = \gy_i, \q\text{ for } i\in\ot.
\end{equation*}

If $L_0$, $L_1$, $L_2$ have zero and both $\gh_1$, $\gh_2$ preserve
the zero, then $L$ can be chosen to have a zero and $\gf_1$, $\gf_2$ can be
chosen to preserve the zero.

If $L_1$ and $L_2$ are finite, then $L$ can be chosen to be finite and
atomistic.
\end{theorem}

This theorem is an extension of T\r uma's theorem---we need only observe
that if the $\gy_i$ are injective, then the $\gf_i$ must be lattice
embeddings. This fact follows from the elementary fact that a lattice
homomorphism
$\gf\colon K\to L$ is an embedding if{f} $\Con{\gf}$ \emph{separates zero},
that is, if{f}
\begin{equation*}
(\Con{\gf})\gQ = \go_L\q\text{implies that}\q\gQ=\go_K,
\q\text{for all }\gQ \in \Con K.
\end{equation*}

We shall apply Theorem~\ref{T:Tumagen} to prove the following strong 
form of Huhn's
Theorem:

\begin{theorem}\label{T:aleph1}
Let $D$ be a distributive algebraic lattice. If $D$ has at most
$\aleph_1$ compact elements, then there exists a locally finite, relatively
complemented lattice $L$ with zero such that $\Con L \iso D$.
\end{theorem}

A lattice $L$ is \emph{congruence-finite}, if $\Con L$ is finite; it is
\emph{$\go$-congruence-finite}, if $L$ can be written as a union,
  \[
  L=\UUm{L_n}{n<\go},
  \]
where $\famm{L_n}{n<\go}$ is an increasing sequence of congruence-finite
sublattices of $L$.

We also apply Theorem~\ref{T:Tumagen} to prove the following:

\begin{theorem}\label{T:EmbLCF}
Every $\go$-congruence-finite lattice $K$ has a
$\go$-congruence-finite, relatively complemented congruence-preserving
extension $L$. Furthermore, if $K$ has a zero, then $L$ can be taken to have
the same zero.
\end{theorem}

\section{Preliminaries}\label{S:SimpleCoat}

\subsection{Notation}
Let $\mathbf{M}_3$ be the five-element nondistributive
modular lattice and let $\mathbf{2}$ denote the two-element chain.

For any lattice $K$, we denote the set of join-irreducible
elements of $K$ by $\J(K)$.

\subsection{Sectionally complemented lattices}
We start with the following stronger form of Dilworth' Theorem (G. 
Gr\"atzer and E.T. Schmidt  \cite{GS98}).

  \begin{theorem}\label{T:seccomp}
Let $D$ be a finite distributive lattice.  Then there exists a finite
sectionally complemented lattice $L$ such that $\Con L$ is isomorphic 
to~$D$ under
an isomorphism $\ga$. Moreover, $L$ contains a Boolean ideal 
generated by the atoms
  \begin{equation*}
\setm{d_p}{p\in \J(D)},
\end{equation*}
and under $\ga$, the congruence $\gQ_L(d_p,0)$ corresponds to $p$, 
for each $p \in
\J(D)$.
  \end{theorem}

Let $K$ and $L$ be lattices, let $f\colon K\to L$ be a lattice
homomorphism. We say that $f$ is \emph{relatively complemented}, if for all
$a$, $b$, $c\in K$ such that $a\leq b\leq c$, there exists a relative
complement of $f(b)$ in the interval $[f(a),f(c)]$ of $L$.

If $f$ is the inclusion map from a sublattice $K$ to the lattice $L$, we say
that $K$ is \emph{relatively complemented} in $L$.

We need the following embedding results:

  \begin{lemma}\label{L:seccomp}\hfill
\begin{enumerate}
    \item For every lattice $L$, there is a bounded, simple, 
sectionally complemented
extension $S(L)$ of $L$ with a dual atom $p$ such that $L$ is relatively
complemented in $S(L)$.
\item If $L$ is finite, then there is a finite, simple,
sectionally complemented extension $S(L)$ of $L$ with a dual atom $p$ 
such that $L$
is relatively complemented in $S(L)$.
\end{enumerate}
  \end{lemma}

For general lattices, by P.M. Whitman \cite{pW46}, every lattice can 
be embedded
in a partition lattice and by O. Ore \cite{oO42}, a partition lattice is
simple and sectionally complemented; it also obviously has a dual 
atom. The second
statement of Lemma~\ref{L:seccomp} follows from
the following very deep result of P.
Pudl\'ak and J. T\r uma \cite{PuTu}: \emph{Every finite lattice can be embedded
into a finite partition lattice}.

  In G. Gr\"atzer and E.T. Schmidt \cite{GS98}, it is pointed out that a
version of this lemma can be proved almost trivially. The comment, and the
simpler proof in \cite{GS98}, also applies to the present version.

\subsection{Congruence-preserving extension}
  Let $L$ be a finite lattice.  A finite lattice $K$ is a
\emph{congruence-preserving extension} of $L$, if $K$ is an extension and every
congruence of $L$ has \emph{exactly one} extension to $K$.  Of course, then the
congruence lattice of $L$ is isomorphic to the congruence lattice of $K$.

A major research tool was discovered by M. Tischendorf \cite{Tisch}:

  \begin{tisch}
  Every finite lattice has a congruence-pre\-serv\-ing extension to a finite
\emph{atomistic} lattice.
  \end{tisch}

A much stronger result was proved in G. Gr\"atzer and E.T. Schmidt 
\cite{GS98}:

  \begin{theorem}\label{T:congpres}\hfill
\begin{enumerate}
\item Every finite lattice has a congruence-pre\-serv\-ing extension 
to a finite
\emph{sectionally complemented}
lattice.
\item Every congruence-finite lattice has a congruence-pre\-serv\-ing 
extension to a
\emph{sectionally complemented}
lattice.
\end{enumerate}
  \end{theorem}

  In the first statement, we cannot strengthen ``sectionally complemented'' to
``relatively complemented'', because the congruence lattice of a 
finite relatively
complemented lattice is always Boolean.

\subsection{$k$-ladders}\label{S:ladder}
Let $k$ be a positive integer. A \emph{$k$-ladder} is a lattice $L$ 
such that, for
any $a \in L$,
\begin{enumerate}
\item $(a]$ is finite;
\item $a$ covers at most $k$ elements.
\end{enumerate}

Note that every $k$-ladder has breadth at most $k$ (see, for example, G.
Gr\"atzer \cite{GLT2} for the definition of breadth).

Every finite chain is a $1$-ladder. The chain $\go$ of all non-negative
integers is also a $1$-ladder.
Note that $k$-ladders are called \emph{$k$-frames} in H.~Dobbertin 
\cite{Dobb86}.

By using the Kuratowski Free Set Theorem, see \cite{Kura}, one can easily
prove that every $k$-ladder has at most $\aleph_{k-1}$ elements, see
S.Z. Ditor \cite{Dito84}. See also H. Dobbertin \cite{Dobb86} for 
the case $k=2$
(his proof does not use the Kuratowski Free Set Theorem). The converse is
obviously true for $k=1$; also for $k=2$, by the following result of 
S.Z. Ditor \cite{Dito84} and by H. Dobbertin
\cite{Dobb86}:

\begin{proposition}\label{P:2ladd}
There exists a $2$-ladder of cardinality $\aleph_1$.
\end{proposition}

\begin{proof}
For $\gx<\go_1$ (the first uncountable ordinal), we construct
inductively the lattices $L_\gx$ with no largest element, as follows. 
Put $L_0=\go$.
If $\gl$ is  countable limit ordinal, put 
$L_\gl=\UUm{L_\gx}{\gx<\gl}$. So assume
that we have constructed $L_\gx$, a countable $2$-ladder with no largest
element. Then $L_\gx$ has a strictly increasing, countable, cofinal, sequence
$\famm{a_n}{n<\go}$. Let $\famm{b_n}{n<\go}$ be a strictly increasing
countable chain, with $b_n\nin L_\gx$, for all $n$. Define $L_{\gx+1}$ by
  \[
  L_{\gx+1}=L_\gx\uu\setm{b_n}{n<\go},
  \]
endowed with the least partial ordering containing the ordering of
$L_\gx$, the natural ordering of $\setm{b_n}{n<\go}$, and all pairs
$a_n<b_n$, for $n<\go$. It is easy to verify that
$L=\UUm{L_\gx}{\gx<\go_1}$ is a $2$-ladder of cardinality $\aleph_1$.
\end{proof}

\section{Proving Theorem~\ref{T:Tumagen}} \label{S:proving}

We prove Theorem~\ref{T:Tumagen} in several steps.

\subsection{Theorem~\ref{T:Tumagen} for $D = \mathbf{2}$} \label{S:two}
In this section, let $D = \mathbf{2}$.

We first state and prove the following special case:

\begin{lemma} \label{L:spec}
Let $L'_0$, $L'_1$, $L'_2$ be lattices and let $\gh'_1\colon L'_0\to L'_1$ and
$\gh'_2\colon L'_0\to L'_2$ be lattice embeddings. Let $D$ be the two-element
chain, and, for $i\in\ot$, let $\gy'_i\colon \Con{L_i}\to D$
satisfy
\begin{equation*}
\gy'_i\gQ = 0_D \q\text{if{f}}\q\gQ=\go_{L'_i}.
\end{equation*}
There is then a lattice $L$ with $1$ and with a dual
atom, there are lattice embeddings
$\gf'_i\colon L'_i\to L$, for $i\in\ot$, with
\begin{equation*}
\gf'_1\compose\gh'_1=\gf'_2\compose\gh'_2,
\end{equation*}
and there is an isomorphism $\ga\colon\Con{L}\to D$ such that
\begin{equation*}
\ga\compose\Con{\gf'_i} = \gy'_i \q\text{for } i\in\ot.
\end{equation*}

If $L'_0$, $L'_1$, $L'_2$ have zero and both $\gh'_1$, $\gh'_2$ preserve
the zero, then $L$ can be chosen to have a zero and $\gf'_1$, $\gf'_2$ can be
chosen to preserve the zero.

If $L'_1$ and $L'_2$ are finite, then $L$ can be chosen to be finite.
\end{lemma}

\begin{proof}
There is a lattice $K$ amalgamating $L'_1$, $L'_2$ over $L'_0$. If $L'_0$,
$L'_1$, $L'_2$ have zero and $\gh'_1$, $\gh'_2$ preserve the zero, then
we can choose $K$ so that $L'_1$ and $L'_2$ are zero-preserving sublattices
of $K$. Observe, also, that if $L'_1$ and $L'_2$ are finite, then $K$
can be chosen finite.

As we pointed out in Lemma~\ref{L:seccomp},
we can embed $K$ into a simple lattice
$L$ that has a 1 and a dual atom, where this embedding preserves the zero,
if $K$ has a zero, and where $L$ is finite, if $K$ is.

For each $i\in\ot$, let $\gf'_i\colon L'_i\to L$ be the composition of
the embedding of $L'_i$ into $K$ with the embedding of $K$ into $L$.
Then
\begin{equation*}
\gf'_1\compose\gh'_1=\gf'_2\compose\gh'_2.
\end{equation*}
Since $L$ is simple, we have an isomorphism $\ga\colon\Con{L}\to D$
such that
\begin{equation*}
\ga\gQ=0_D\q\text{if{f}}\q\gQ=\go_L.
\end{equation*}
For each $i\in\ot$ and each $\gQ\in\Con{L'_i}$,
\begin{equation*}
(\Con{\gf'_i})\gQ = \go_L\q\text{if{f}}\q\gQ=\go_{L'_i},
\end{equation*}
since $\gf'_i$ is an embedding.

Thus,
\begin{equation*}
\ga\compose\Con{\gf'_i} = \gy'_i,
\end{equation*}
concluding the proof of the lemma.
\end{proof}

We proceed to prove Theorem~\ref{T:Tumagen} for $D = \mathbf{2}$. For 
each $i\in\ot$,
set \begin{equation*}
\gQ_i = \JJ\famm{\gQ\in \Con{L_i}}{\gy_i\gQ = 0_D},
\end{equation*}
and set
\begin{equation*}
\gQ_0=\setm{\vv<x,y>\in L_0}{\gh_1x\equiv{\gh_1y}\,(\gQ_1)} =
\setm{\vv<x,y>\in L_0}{{\gh_2x}\equiv{\gh_2y}\,(\gQ_2)}.
\end{equation*}
For each
$i\in\set{0,1,2}$, set $L'_i=L_i/\gQ_i$ and let
$\gp_i\colon L_i\twoheadrightarrow L'_i$
be the canonical surjection. Note that $\gQ_i = \ker{\gp_i}$. We then
have lattice embeddings $\gh'_1\colon L'_0\to L'_1$ and
$\gh'_2\colon L'_0\to L'_2$ such that
\begin{equation*}
\gp_i\compose\gh_i = \gh'_i\compose\gp_0, \q \text{ for } i\in\ot.
\end{equation*}
Furthermore, we have mappings
$\gy'_1\colon \Con{L'_1}\to D$ and $\gy'_2\colon \Con{L'_2}\to D$
with
\begin{equation*}
\gy'_i\gQ = 0_D \q \text{ if{f} }\q \gQ=\go_{L'_i}, \q \text{ for } i\in\ot
\end{equation*}
such that
\begin{equation*}
\gy'_i\compose \Con{\gp_i} = \gy_i, \q \text{ for } i\in\ot.
\end{equation*}

We apply Lemma~\ref{L:spec} to get the lattice $L$, the embeddings
$\gf'_i\colon L'_i\to L$, and the isomorphism $\ga\colon\Con{L}\to D$
with
\begin{equation*}
\ga\compose\Con{\gf'_i} = \gy'_i, \q\text{ for } i\in\ot.
\end{equation*}
For each $i\in\ot$, set
\begin{equation*}
\gf_i = \gf'_i\compose\gp_i\colon L_i\to L.
\end{equation*}
Then
\begin{align*}
\gf_1\compose\gh_1 &= \gf'_1\compose\gp_1\compose\gh_1 =
\gf'_1\compose\gh'_1\compose\gp_0\\ &= \gf'_2\compose\gh'_2\compose\gp_0 =
\gf'_2\compose\gp_2\compose\gh_2 = \gf_2\compose\gh_2,
\end{align*}
and
\begin{equation*}
\ga\compose\Con{\gf_i}=\ga\compose(\Con{\gf'_i})\compose(\Con{\gp_i})
=\gy'_i\compose\Con{\gp_i}=\gy_i,
\end{equation*}
for $i\in\ot$, concluding the proof of Theorem~\ref{T:Tumagen} for $D 
= \mathbf{2}$.

\subsection{Theorem~\ref{T:Tumagen} for $D$ Boolean}

In this section, let $D$ be a finite Boolean lattice. We prove
Theorem~\ref{T:Tumagen} with the following addition:

\begin{add}
$L$ contains a Boolean dual ideal isomorphic to $D$ with a set
\begin{equation*}
\setm{d_p}{p\in \J(D)}
\end{equation*}
its set of dual atoms. For each $p\in \J(D)$,
\begin{equation*}
\ga\gQ_L(d_p,1)=p.
\end{equation*}
\end{add}

\begin{proof}
The set $\J(D)$ is the set of atoms of $D$.
For each $p\in \J(D)$,
we have a zero-preserving lattice surjection $\gb_p\colon D \to \two$
such that $\gb_p(x)=1$
if{f} $p\leq x$. Then
\begin{equation*}
\gb = \prod\famm{\gb_p}{p\in \J(D)}\colon D\to\prod\famm{\two}{p\in \J(D)}
\end{equation*}
is an isomorphism.

For each $p\in \J(D)$, set $\gy_{pi}=\gb_p\compose\gy_i$, for $i\in\ot$ and
apply the case $D = \mathbf{2}$ to the configuration $\gh_i\colon L_0\to L_i$,
$\gy_{pi}\colon\Con{L_i}\to\two$ to obtain a simple lattice $L_p$ 
with a $1$ and
a dual atom $d'_p$, lattice homomorphisms $\gf_{pi}\colon L_i\to L_p$
with $\gf_{p1}\compose\gh_1=\gf_{p2}\compose\gh_2$, and an isomorphism
$\ga_p\colon\Con{L_p}\to\two$ with
\begin{equation*}
\ga_p\compose\Con{\gf_{pi}}=\gy_{pi}=\gb_p\compose\gy_i.
\end{equation*}

We then set
\begin{equation*}
L=\prod\famm{L_p}{p\in \J(D)}
\end{equation*}
and set
\begin{equation*}
\gf_i = \prod\famm{\gf_{pi}}{p\in \J(D)} \colon L_i\to L.
\end{equation*}
Then
\begin{equation*}
\gf_1\compose\gh_1=\gf_2\compose\gh_2.
\end{equation*}
Now,
\begin{equation*}
\Con{\gf_i}=\prod\famm{\Con{\gf_{pi}}}{p\in \J(D)}.
\end{equation*}
Thus,
\begin{align*}
\prod\famm{\ga_p}{p\in \J(D)}\compose\Con{\gf_i}
&= \prod\famm{\ga_p\compose\Con{\gf_{pi}}}{p\in \J(D)} \\
&= \prod\famm{\gb_p\compose\gy_i}{p\in \J(D)} \\
&= \prod\famm{\gb_p}{p\in \J(D)}\compose\gy_i  \\
&= \gb\compose\gy_i.
\end{align*}
Setting $\ga=\gb^{-1}\compose\prod\famm{\ga_p}{p\in \J(D)}$, we thus get
an isomorphism $\ga\colon\Con{L}\to D$ with
\begin{equation*}
\ga\compose\Con{\gf_i} = \gy_i.
\end{equation*}

For each $q\in \J(D)$, we define $d_q\in L =\prod\famm{L_p}{p\in \J(D)}$ by
setting
\begin{equation*}
(d_q)_p =
\begin{cases}
d'_q, &\text{if $p=q$;} \\
1_{L_p}, &\text{otherwise.}
\end{cases}
\end{equation*}
Then each $d_q$ is a dual atom of $L$, and the dual ideal of $L$ generated
by $\setm{d_p}{p\in \J(D)}$ is
\begin{equation*}
\prod\famm{\set{d'_p,1_{L_p}}}{p\in \J(D)},
\end{equation*}
a Boolean lattice with $\setm{d_p}{p\in \J(D)}$ its set of dual atoms.

Now, $\Con{L}=\prod\famm{\Con{L_p}}{p\in \J(D)}$ and each $L_p$ is simple.
Thus, for $p$, $q\in \J(D)$,
the $p$-th component of $\gQ_L(d_q,1_L)$ satisfies
\begin{equation*}
(\gQ_L(d_q,1_L))_p =
\begin{cases}
\gQ_{L_q}(d'_q,1_{L_q}) = \gi_{L_p}, &\text{if $p=q$;} \\
\gQ_{L_p}(1_{L_p},1_{L_p}) = \go_{L_p}, &\text{otherwise.}
\end{cases}
\end{equation*}
Then, for each $p\in \J(D)$,
\begin{equation*}
\gb_pq = \ga_p(\gQ_L(d_q,1))_p,
\end{equation*}
that is,
\begin{equation*}
\gb q = \prod\famm{\ga_p}{p\in \J(D)}\gQ_L(d_q,1),
\end{equation*}
that is,
  \begin{equation*}
  q = \ga\gQ_L(d_q,1).
  \end{equation*}

Since finite direct products preserve the zero and finiteness, the proof is
completed.
\end{proof}

\subsection{The general proof}\label{S:main}

We let $B$ be the Boolean lattice generated by $D$, and let
$\gh\colon D\hookrightarrow B$ be the canonical embedding.
For each $x\in B$, let $\gr x$ denote
the smallest element of $D$ containing $x$. Then we get a \jz-homomorphism
$\gr\colon B\to D$ such that
  \begin{equation}\label{E:retract}
  \gr\compose\gh = \id_D.
  \end{equation}

Note that $\gr\res_{\J(B)}$ is just the usual
dual of $\gh$ in the duality between distributive lattices and posets.
In our case of $B$ being the Boolean lattice generated by $D$, we get
an isotone bijection
\begin{equation} \label{E:bij}
\gr\res_{\J(B)}\colon \J(B)\to \J(D).
\end{equation}

We apply the special case where
$D$ is Boolean to the system $L_0$, $L_1$, $L_2$, $B$
with the complete $\JJ$-homomorphisms $\gh\compose\gy_i
\colon\Con{L_i}\to B$, $i\in\ot$, and obtain a lattice $K_0$ and
lattice homomorphisms $\gf'_i\colon L_i\to K_0$, $i\in\ot$, satisfying
\begin{equation*}
\gf'_1\compose\gh_1 = \gf'_2\compose\gh_2,
\end{equation*}
and an isomorphism
\begin{equation*}
\ga_0\colon\Con{K_0}\to B
\end{equation*}
such that
\begin{equation} \label{E:psi}
\ga_0\compose\Con{\gf'_i} = \gh\compose\gy_i,\q\text{for } i\in\ot.
\end{equation}
Furthermore, $K_0$ contains a finite Boolean dual ideal $H$ with $|\J(B)|$ dual
atoms $d'_p$, $p\in \J(B)$, such that
\begin{equation*}
\ga_0\gQ_{K_0}(d'_p,1) = p,
\end{equation*}
for each $p\in \J(B)$.

By Theorem \ref{T:seccomp}, there is a finite lattice $K_1$ and there is an
isomorphism
\begin{equation*}
\ga_1\colon\Con{K_1}\to D
\end{equation*}
such that $K_1$ contains a Boolean ideal $I$ with $|\J(D)|$ dual atoms $d_p$,
$p\in\J(D)$, and
\begin{equation*}
\ga_1\gQ_{K_1}(d_p,1_I) = p,
\end{equation*}
for each $p \in \J(D)$.

In view of the bijection \eqref{E:bij}, there is an isomorphism of the
dual ideal $H$ of $K_0$ with the ideal $I$ of $K_1$, whereby $d'_p\in H$
corresponds to $d_{\gr p}$, for each $p\in \J(B)$. We let $L$
be the lattice obtained by gluing $K_1$ to the top of $K_0$ by identifying
$H$ with $I$ under this isomorphism, so that $K_1$ is a subset of $L$.
We then have an embedding
$\ge_0\colon K_0\to L$, where
\begin{equation*}
\ge_0\colon d'_p\mapsto d_{\gr p},
\end{equation*}
for $p\in \J(B)$, and an embedding $\ge_1\colon K_1\to L$, where
\begin{equation*}
\ge_1\colon d_p\mapsto d_p,
\end{equation*}
for $p\in \J(D)$. Then $\Con{\ge_1}\colon \Con{K_1}\to\Con{L}$, whereby
\begin{equation} \label{E:cone1}
\Con{\ge_1}\colon\gQ_{K_1}(d_p,1_I)\mapsto \gQ_L(d_p,1_I),
\end{equation}
is an isomorphism, and the \jz-homomorphism
$\Con{\ge_0}\colon \Con{K_0}\to\Con{L}$ satisfies
\begin{equation} \label{E:cone0}
\Con{\ge_0}\colon \gQ_{K_0}(d'_p, 1_{K_0})\mapsto \gQ_L(d_{\gr p}, 1_I).
\end{equation}

For each $i\in\ot$, we set
\begin{equation*}
\gf_i = \ge_0\compose\gf'_i\colon L_i\to L.
\end{equation*}
Then
\begin{equation*}
\gf_1\compose\gh_1 = \ge_0\compose\gf'_1\compose\gh_1 =
\ge_0\compose\gf'_2\compose\gh_2 = \gf_2\compose\gh_2.
\end{equation*}

We have an isomorphism $\ga\colon\Con{L}\to D$ defined by
\begin{equation*}
\ga = \ga_1\compose(\Con{\ge_1})^{-1}.
\end{equation*}
We proceed to show that $\ga\compose\Con{\gf_i} = \gy_i$, for $i\in\ot$.
By the definition of $\ga$ and $\gf_i$,
\begin{equation} \label{E:alpha}
\ga\compose\Con{\gf_i} = \ga_1\compose(\Con{\ge_1})^{-1}
\compose(\Con{\ge_0})\compose(\Con{\gf'_i}).
\end{equation}
By \eqref{E:cone1} and \eqref{E:cone0},
\begin{equation*}
(\Con{\ge_1})^{-1}\compose(\Con{\ge_0})\colon \gQ_{K_0}(d'_p,1_{K_0})
\mapsto \gQ_{K_1}(d_{\gr p},1_I),
\end{equation*}
for each $p\in \J(B)$. Thus,
\begin{equation*}
\ga_1\compose(\Con{\ge_1})^{-1}\compose(\Con{\ge_0})\compose\ga_0^{-1}
\colon p\mapsto\gr p,
\end{equation*}
for each $p\in \J(B)$. Therefore,
  \begin{equation} \label{E:rhoeq}
  \ga_1\compose(\Con{\ge_1})^{-1}\compose(\Con{\ge_0})\compose\ga_0^{-1}
  = \gr,
  \end{equation}
since both sides are \jz-homomorphisms. Thus, for $i\in\ot$,
\begin{align*}
\ga\compose\Con{\gf_i} &= \gr\compose\ga_0\compose(\Con{\gf'_i}),
&\text{by \eqref{E:alpha} and \eqref{E:rhoeq},} \\
&= \gr\compose\gh\compose\gy_i, &\text{by \eqref{E:psi},} \\
&= \gy_i, &\text{by \eqref{E:retract}.}
\end{align*}

This concludes the proof for arbitrary lattices $L_i$ and homomorphisms
$\gh_i$.

We note that if $K_0$ has a zero, then $\ge_0\colon K_0\to L$
preserves the zero. Thus, by the special case $D$ is Boolean, if the 
$L_i$ each have
a zero and if the $\gh_i$ preserve the zero, then the $\gf'_i$ and, 
consequently,
the $\gf_i$ preserve the zero.

We note, also, that if $L_1$ and $L_2$ are finite, then so is $K_0$ and
thus so is $L$. Then, using Tischendorf's Theorem, we can
replace $L$ by a finite atomistic lattice.

This concludes the proof of the Theorem~\ref{T:Tumagen}.

\section{Proving Theorem \ref{T:aleph1}}

\subsection{Congruence-preserving extensions}
We shall now establish two results, the first a strengthening of both
parts of Theorem~\ref{T:congpres}, and the second a strengthening of
Theorem~\ref{T:congpres}(ii) in a different direction:

\begin{lemma}\label{L:HLadded}
Let $K$ be a congruence-finite lattice. Then $K$ has a
congruence-preserving relatively complemented embedding into a sectionally
complemented lattice $K'$. If $K$ has a zero, then one can assume that
$K'$ has the same zero. If $K$ is finite, then $K'$ can be chosen to be
finite.
\end{lemma}

\begin{proof}[Outline of proof]
We follow the original proof in G. Gr\"atzer and E.T. Schmidt 
\cite{GS98}, with just one small addition.  If $K$
is a congruence-finite lattice, the congruence-preserving
sectionally complemented extension of $K$ is constructed as follows. 
Since $\Con K$
is a finite distributive lattice, we can associate with it the finite 
sectionally
complemented lattice $L_0$ of Theorem
\ref{T:seccomp} such that $\Con L_0\iso\Con K$. On the
other hand, denote by $\M(\Con K)$ the set of all meet-irreducible 
congruences of
$K$. The \emph{rectangular extension} of $K$ is defined by
  \[
  \RR(K)=\prod\famm{K/\gQ}{\gQ\in\M(\Con K)}.
  \]
Let $K_\gQ$ be a simple sectionally complemented extension of $K/\gQ$
such that, in addition,  $K/\gQ$ is relatively complemented in $K_\gQ$
(we use Lemma~\ref{L:seccomp}). If $K$ is finite we choose $K_\gQ$ finite.
Put
  \[
  \RRh(K)=\prod\famm{K_\gQ}{\gQ\in\M(\Con K)}.
  \]
Note that the diagonal map from $K$ into $\RRh(K)$, that sends every $x\in K$
to $\famm{[x]\gQ}{\gQ\in\M(\Con K)}$, has the congruence extension 
property, but
it is not necessarily congruence-preserving (the congruence lattice of
$\RRh(K)$ is Boolean). However, the sectionally complemented extension $K'$
constructed in \cite{GS98} is obtained by considering the lattice of finitely
generated ideals of the chopped lattice $L_0\uu\RRh(K)$, with the two
isomorphic Boolean
sublattices of $L_0$ and $\RRh(K)$ identified. Since $K$ is already relatively
complemented in $\RRh(K)$, it is \emph{a fortiori} relatively complemented in
$K'$.

If $K$ has a zero, then the above construction preserves this zero.
Furthermore, if $K$ is finite, then $K'$ is finite. \end{proof}

\begin{theorem}\label{T:GSimpr}
Let $K$ be a congruence-finite lattice.
Then $K$ has a
congruence-preserving embedding into a relatively complemented lattice
$L$. Furthermore, if $K$ has a zero, then one can assume that $L$ has the same
zero.
\end{theorem}
\begin{proof}
We use Lemma~\ref{L:HLadded} to construct a sequence
$\famm{K^{(n)}}{n<\go}$ of lattices. Set $K^{(0)}=K$, and, proceeding
inductively, for each $n$ set $K^{(n+1)}=(K^{(n)})'$, the lattice
$(K^{(n)})'$ being the extension of $K^{(n)}$ guaranteed by
Lemma~\ref{L:HLadded}.
To conclude the proof, it suffices to take
$L=\UUm{K^{(n)}}{n<\go}$. \end{proof} 

\subsection{Proving Theorem \ref{T:aleph1}}
\label{S:aleph1}
Let $S$ be the \jz-semilattice of all compact elements of $D$.
By definition, $S$ is distributive.
By P. Pudl\'ak's Lemma (see \cite{Pudl85}),
every finite subset of $S$
is contained in a finite distributive \jz-subsemilattice of $S$. We use this to
construct a direct system of finite distributive subsemilattices of $S$ as
follows. First, by Proposition~\ref{P:2ladd}, there exists a $2$-ladder of
cardinality
$\aleph_1$, say, $\vv<I,\leq>$. Let $\gp\colon I\twoheadrightarrow S$ be a
surjective map such that $\gp(0_I)=0_S$.
We define a family $\famm{S_i}{i\in I}$ of finite distributive
\jz-subsemilattices of
$S$, as follows. We put $S_0=\set{0}$, and, for all $i\in I$, we let $S_i$ be
a finite distributive \jz-subsemilattice of $S$ containing the subset
  \[
  \UUm{S_j}{j<i}\uu\set{\gp(i)}.
  \]
Since $\gp(0_I)=0_S$, we can take $S_0=\set{0}$.
Then $S$ is the directed union of all $S_i$, for $i\in I$.
We denote by $\gf^i_j$ the inclusion map from $S_i$
into $S_j$, for all $i\leq j$ in $I$.

Let $\gr\colon I\to\go$ be any strictly increasing map from $I$ to $\go$ (for
example, the height function on $I$). We put
  \[
  I_n=\setm{i\in I}{\gr(i)\leq n},
  \]
for all $n<\go$. By induction on $n$, we construct a family of finite lattices
$L_i$, maps $\ge_i\colon\Con L_i\to S_i$, for $i\in I_n$, and
$\set{0}$-lattice homomorphisms $f^i_j\colon L_i\to L_j$, for $i\leq j$ in
$I_n$, satisfying the following properties:

\begin{itemize}
\item[(a)] $f^i_i=\id_{L_i}$, for all $i\in I_n$.

\item[(b)] $f^i_k=f^j_k\circ f^i_j$, for all $i$, $j$, $k\in I_n$.

\item[(c)] $\ge_i$ is an isomorphism from $\Con L_i$ onto $S_i$, for all
$i\in I_n$.

\item[(d)] $\ge_j\circ\Con f^i_j=\gf^i_j\circ\ge_i$, for all $i\leq j$ in
$I_n$.

\item[(e)] $f^i_j[L_i]$ is relatively complemented in $L_j$, for all $i<j$ in
$I_n$.
\end{itemize}

For $n=0$, we just take $L_0=\set{0}$ (because $S_0=\set{0}$). Let us 
assume that
we have performed the construction at level $n$; we show how to 
extend it to the
level $n+1$. So, let $i\in I_{n+1}-I_n$. Since $I$ is a $2$-ladder, $i$ has (at
most) two immediate predecessors in $I$, say, $i_1$ and $i_2$. Note that $i_1$
and $i_2$ need not be distinct. For $k\in\ot$, the map
  \[
  \gy_k=\gf^{i_k}_i\circ\ge_{i_k}
  \]
is a \jz-embedding from $\Con L_{i_k}$ into $S_i$, and the equality
  \[
  \gy_1\circ\Con f^{i_1\mm i_2}_{i_1}=\gy_2\circ\Con f^{i_1\mm i_2}_{i_2}
  \]
holds. By Theorem~\ref{T:Tumagen}, there is a finite lattice $L_i$, there are
$\set{0}$-lattice homomorphisms $g_k\colon L_{i_k}\to L_i$, for $k\in\ot$, and
and there is an isomorphism $\ge_i\colon\Con L_i\to S_i$ such that
  \begin{gather}
  g_1\circ f^{i_1\mm i_2}_{i_1}=g_2\circ f^{i_1\mm i_2}_{i_2},\label{Eq:g12f}\\
  \ge_i\circ\Con g_k=\gy_k,\qquad\text{for }k\in\ot,\label{Eq:geigk}
  \end{gather}
hold.
Furthermore, if $i_1=i_2$, then replacing $g_2$ by $g_1$ does not change the
validity of \eqref{Eq:g12f} and \eqref{Eq:geigk}. Thus we may define
$f^{i_k}_i=g_k$, for
$k\in\ot$, and \eqref{Eq:g12f}, \eqref{Eq:geigk} take the following form:
  \begin{gather}
  f^{i_1}_i\circ f^{i_1\mm i_2}_{i_1}=f^{i_2}_i\circ f^{i_1\mm i_2}_{i_2},
  \label{Eq:newf12f}\\
  \ge_i\circ\Con f^{i_k}_i=\gy_k,\qquad\text{for }k\in\ot.\label{Eq:newfs}
  \end{gather}
Furthermore, by one application of Lemma~\ref{L:HLadded}, we may
assume that both embeddings $f^{i_1}_i$ and $f^{i_2}_i$ are relatively
complemented (in $L_i$). This takes care of (e).

So we have defined $f^j_i$, for all $i$ and $j$ in $I_{n+1}$ such that $j$ is
an immediate predecessor of $i$ in $I_{n+1}$. We extend this definition to
arbitrary $i$, $j\in I_{n+1}$ such that $j\leq i$. If $j=i$, then we put
$f^i_j=\id_{L_i}$. Now assume that $j<i$ in $I_{n+1}$, with $i\nin I_n$.
There exists an index $k\in\ot$ such that $j\leq i_k$.
The only possible choice for $f^j_i$ is to define it as
  \begin{equation}\label{Eq:extf}
  f^j_i=f^{i_k}_i\circ f^j_{i_k},
  \end{equation}
except that this should be independent of $k$. This means that if
$j\leq i_1\mm i_2$, then the equality
  \begin{equation}\label{Eq:cohfij}
  f^{i_1}_i\circ f^j_{i_1}=f^{i_2}_i\circ f^j_{i_2}
  \end{equation}
should hold. We compute:
  \begin{align*}
  f^{i_1}_i\circ f^j_{i_1}&=
  f^{i_1}_i\circ f^{i_1\mm i_2}_{i_1}\circ f^j_{i_1\mm i_2}\\
  &=f^{i_2}_i\circ f^{i_1\mm i_2}_{i_2}\circ f^j_{i_1\mm i_2}\text{\qq \qq (by
\eqref{Eq:newf12f})}\\
  &=f^{i_2}_i\circ f^j_{i_2},
  \end{align*}
which establishes \eqref{Eq:cohfij}.

At this point, the $\set{0}$-lattice embeddings $f^j_i\colon L_j\to L_i$ are
defined for all $j\leq i$ in $I_{n+1}$. The verification of conditions
(a)--(c) above is then straightforward. Let us verify (d). Let $i\leq j$ in
$I_{n+1}$, we prove that
\begin{equation}\label{Eq:fijeps}
\gf^i_j\circ\ge_i=\ge_j\circ\Con f^i_j.
\end{equation}
The only nontrivial case happens if $j\in I_{n+1}-I_n$ and $i<j$. It
suffices then to verify \eqref{Eq:fijeps} for the pairs $\vv<i,j_*>$ and
$\vv<j_*,j>$, where $j_*$ is any immediate predecessor of $j$ such that
$i\leq j_*$. For the pair $\vv<i,j_*>$, this follows
from the induction hypothesis, while for the pair $\vv<j_*,j>$, this follows
from \eqref{Eq:newfs}.

Hence the construction of the $L_i$, $\ge_i$, $f^i_j$ is carried out for
the whole poset $I$. Let $L$ be the direct limit of all the $L_i$, $i 
\in I$, with
the transition maps $f^i_j$, for $i\leq j$ in $I$. Then $\Conc L$ is the direct
limit of the $\Conc L_i$, with the transition maps $\Conc f^i_j$, in the
category of distributive\jz-semilattices and \jz-homomorphisms. Thus,
by (c) and (d), $\Conc L$ is isomorphic to the direct limit of the $S_i$ with
the transition maps $\gf^i_j$, for $i\leq j$ in $I$. Hence, $\Conc L\iso S$,
from which it follows that $\Con L\iso D$. The fact that $L$ is relatively
complemented follows from condition (e) above.

\section{Proving Theorem \ref{T:EmbLCF}}

By definition, $K$ can be written as a union,
  \[
  K=\UUm{K_n}{n<\go},
  \]
where $\famm{K_n}{n<\go}$ is an increasing sequence of congruence-finite
sublattices of $K$. Furthermore, if $K$ has a zero, then we can assume that
$0$ belongs to $K_n$, for all $n<\go$. Denote by $e_n$ the inclusion map from
$K_n$ into $K_{n+1}$.
For $n<\go$, let us assume that we have constructed a relatively complemented
lattice $L_n$ and a congruence-preserving embedding
$u_n\colon K_n\hookrightarrow L_n$ such that $u_n$ preserves the zero if $K_n$
has a zero. We apply Theorem~\ref{T:Tumagen} to the lattice homomorphisms
  \[
  u_n\colon K_n\hookrightarrow L_n,\qquad
  e_n\colon K_n\hookrightarrow K_{n+1},
  \]
the semilattice $D=\Con K_{n+1}$,
and the \jz-semilattice homomorphisms
  \begin{align*}
  \gf=\Con u_n&\colon\Con K_n\to\Con L_n,\\
  \gy=(\Con e_n)\circ(\Con u_n)^{-1}&\colon\Con L_n\to\Con K_{n+1}.
  \end{align*}
We obtain a lattice $L_{n+1}$, lattice homomorphisms
  \[
  f_n\colon L_n\to L_{n+1},\qquad
  u_{n+1}\colon K_{n+1}\to L_{n+1},
  \]
and an isomorphism $\ga_n\colon\Con L_{n+1}\to\Con K_{n+1}$ such that
the following equalities hold:
  \begin{align}
  u_{n+1}\circ e_n&=f_n\circ u_n,\label{Eq:isodia}\\
  \ga_n\circ\Con f_n&=(\Con e_n)\circ(\Con u_n)^{-1},\label{Eq:fn0sep}\\
  \ga_n\circ\Con u_{n+1}&=\id_{\Con K_{n+1}}.\label{Eq:un+1CP}
  \end{align}
By Theorem~\ref{T:GSimpr}, one can further assume that
$L_{n+1}$ is relatively complemented.
By \eqref{Eq:un+1CP}, the map $\Con u_{n+1}$ is an isomorphism and so 
$u_{n+1}$ is
congruence-preserving. By \eqref{Eq:fn0sep}, the map $\Con f_n$ separates zero
(because $\Con e_n$ does), that is, $f_n$ is a lattice embedding.

Let $L$ be the direct limit of all the $L_n$, with the transition maps
\[
f_m\circ\cdots\circ f_{n-1}\colon L_m\hookrightarrow L_n,
\]
for $m<n$ in $\go$. Denote by $g_n\colon L_n\to L$ the corresponding limiting
maps. By \eqref{Eq:isodia} and the fact that all the $u_n$ are
congruence-preserving embeddings, the sequence $\famm{u_n}{n<\go}$ defines a
congruence-preserving embedding $u\colon K\to L$ by
  \[
  u(x)=g_n\circ u_n(x),\quad\text{if }x\in K_n,\text{ for }n<\go.
  \]
Since all the $L_n$ are relatively complemented, $L$ is relatively
complemented. If $K$ has a zero, then all the maps $u_n$ and $f_n$
preserve the zero, thus $L$ has the same zero as $K$.

If $K$ is locally finite, then we can assume that all the $K_n$ are finite,
and we can then take all the $L_n$ finite. In particular, $L$ is also locally
finite.  This concludes the proof of Theorem~\ref{T:EmbLCF}.

\section{Discussion}

\subsection{Theorem~\ref{T:Tumagen}}

In Theorems \ref{T:Tumagen}--\ref{T:EmbLCF}, the bound zero is
preserved. We do not know whether the theorems of this paper have analogues for
\emph{bounded} lattices:

\begin{problem}
In the statement of Theorem~\ref{T:Tumagen}, let us assume that 
$L_0$, $L_1$, and
$L_2$ are bounded lattices and that $\gh_1$, $\gh_2$ are 
$\set{0,1}$-preserving.
Can the lattice $L$ of the conclusion be taken bounded, with both $\gf_1$ and
$\gf_2$ $\set{0,1}$-preserving? In addition, if $L_0$, $L_1$, and 
$L_2$ are finite,
can $L$ be taken finite?
\end{problem}

\subsection{Theorem~\ref{T:aleph1}}
 From the results of the third author in \cite{Wehr1,Weurp}, the $\aleph_1$
bound in the statement of Theorem~\ref{T:aleph1} is \emph{optimal}, because
there are algebraic distributive lattices with $\aleph_2$ compact elements
that cannot be represented as congruence lattices of relatively
complemented lattices.

There are stronger forms of Theorem~\ref{T:aleph1}.
For example, a result of K.R. Goodearl and F. Wehrung \cite{GoWe} states that
every distributive \jz-semilattice is the direct limit of a family of finite
\emph{Boolean} \jz-semilattices and \jz-homomorphisms. Since every 
finite lattice
embeds into a finite geometric lattice, one can prove that
\emph{the lattice $L$ of Theorem~\tup{\ref{T:aleph1}} can be assumed 
to be a direct
limit of finite geometric lattices}. Similarly, using
P. Pudl\'ak and J. T\r uma \cite{PuTu}, we can prove that
$L$ can be assumed to be a direct limit of lattices, each of which is a finite
product of finite partition lattices.

In neither of these cases is $L$ \emph{modular}. However, using the results of
\cite{Wehr}, one can show that the lattice $L$ of Theorem~\ref{T:aleph1} can be
taken to be \emph{sectionally complemented and modular}; in addition, in this
case, $L$ can be assumed to be bounded, if the largest element of $D$ 
is compact.
The local finiteness of $L$ is lost.

\begin{problem}
If the lattice $L$ has at most $\aleph_1$ compact congruences, does $L$ have a
relatively complemented congruence-preserving extension.
\end{problem}

A variant of this problem, was raised by the first and the last author
at the August 1998 Szeged meeting:

  \begin{problem}
  Let $L$ be an infinite lattice with $|L|\leq\aleph_1$. Does $L$ have a
  congruence-preserving extension to a (sectionally complemented) relatively
  complemented lattice?
  \end{problem}

\subsection{Theorem~\ref{T:EmbLCF}}

The countability assumption of the statement of Theorem~\ref{T:EmbLCF} is
essential: by M. Plo\v s\v cica, J. T\r uma, and F. Wehrung \cite{PTWe},
the free lattice with $\aleph_2$ generators in the variety generated 
by $\mathbf{M}_3$
(or any finite, nondistributive lattice) does not have a
congruence-preserving embedding into a relatively complemented lattice.

Not every countable lattice is $\go$-congruence-finite: take any finitely
generated, non congruence-finite lattice, for example, the free lattice on $n$
generators, where $n\geq3$.

\begin{problem}

Is it true that every bounded, $\go$-congruence-finite lattice $L$ has a
con\-gruence-pre\-serving extension into a relatively complemented lattice that
preserves the bounds?
\end{problem}

\end{document}